\documentclass[11pt]{article}

\usepackage{booktabs}
\usepackage[T1]{fontenc}
\usepackage{lmodern}
\usepackage[margin=1in]{geometry}
\usepackage{amsmath,amssymb,amsthm,mathtools}
\usepackage{enumitem}
\usepackage{microtype}
\usepackage{xurl}
\usepackage[hidelinks,pdfusetitle,hypertexnames=false]{hyperref}
\hypersetup{pdfauthor={Pedro Poitevin}}

\newtheorem{theorem}{Theorem}[section]
\newtheorem{proposition}[theorem]{Proposition}
\newtheorem{lemma}[theorem]{Lemma}
\newtheorem{corollary}[theorem]{Corollary}
\theoremstyle{definition}
\newtheorem{definition}[theorem]{Definition}

\theoremstyle{remark}
\newtheorem{remark}[theorem]{Remark}
\numberwithin{equation}{section}
\setlist[enumerate]{topsep=3pt,itemsep=0pt,parsep=0pt,partopsep=0pt}
\DeclareMathOperator{\Th}{Th}
\newcommand{\N}{\mathbb N}
\newcommand{\R}{\mathbb R}
\newcommand{\eps}{\varepsilon}
\newcommand{\Ball}{B}

\title{The norm formula in a Banach space and in its complete theory}
\author{Pedro Poitevin\thanks{Mathematics Department, Salem State University,
Salem, MA 01970, USA. Email: \href{mailto:ppoitevin@salemstate.edu}{\texttt{ppoitevin@salemstate.edu}}.}}

\date{November 15, 2026}

\begin{document}
\maketitle

\begin{abstract}
\noindent
Let \(h\) be the Hilbert norm on the real space
\(\ell_2(\N)\oplus_2\ell_2(\N)\), and put \(q(u,v)=\sup_{i<j}|u_i+v_j|\).
For each of Bollob\'as's additive norms \(h+\eps q\), \(\eps>0\), and
for the quadratic norm \(\sqrt{h^2+q^2}\), the norm formula
\(\lVert x+y\rVert\) is unstable on the renormed space and has no
independence property in its complete theory. We approximate \(q(x+y)\)
uniformly on bounded Hilbert balls by truncating small coordinates.
For \(n\) fixed test vectors and thresholds, the number of threshold
patterns of \(q(x+y)\) as \(y\) varies is bounded by a polynomial in \(n\)
when the support sizes of the test vectors and of \(y\) are bounded.
These estimates give finite fat-shattering dimension for the translates
of \(q\) at every positive scale.
The norm formula on a space can be stable, unstable without the independence
property, or have the independence property. We give an example in each
of these three cases where the formula has the independence property in
the complete theory, and an example in each of the first two cases where
it lacks that property in the complete theory. All five spaces satisfy
a common midpoint estimate for uniform convexity, and each admits an
isomorphism to Hilbert space of distortion at most \(\sqrt3\).
\end{abstract}

\section{Introduction}

Can the norm formula \(\varphi(x,y)=\lVert x+y\rVert\) be stable on a
Banach space and have the independence property in its complete theory?
Theorem~\ref{thm:main} gives such a space by a uniformly convex renorming of
\(\ell_2\).
We also ask whether the norm formula for Bollob\'as's unstable renorming
of \(\ell_2\) has the independence property in its complete theory.
Write \(\N=\{1,2,\ldots\}\). On \(H=\ell_2(\N)\oplus_2\ell_2(\N)\),
let \(h\) be the Hilbert norm and put \(q(u,v)=\sup_{i<j}|u_i+v_j|\).
Bollob\'as's additive renorming is \(h+\eps q\) for \(\eps>0\)
\cite{Garling1982,IovinoStableBanachI}. Corollary~\ref{cor:additive-norm}
shows that, for every \(\eps>0\), the norm formula has no independence
property in any model of the complete theory of the renormed space.

Throughout, all Banach spaces are real. We work in continuous first-order
logic for metric structures \cite{BBHU08}. The Banach language is countable
and has a sort for each closed ball of nonnegative rational radius, with
rational scalar multiplication. Each formula is bounded and uniformly
continuous on the product of its variable sorts. We write \(\Ball_X\)
for the closed unit ball of \(X\) and \(R\Ball_X\) for its closed ball of
radius \(R\geq0\).

For a bounded formula \(\varphi(x,y)\) on a metric structure \(M\), write
\(M^x\) and \(M^y\) for the products of the sorts of its two variable tuples.

Following Khanaki \cite[Definition~4.1(i)]{KhanakiBanach}, \(\varphi\) has the
\emph{independence property} (IP) on \(M\) if there is a sequence
\((a_n)_{n\in\N}\) in \(M^x\) and real numbers \(r<s\) such
that, for every pair of finite disjoint sets \(E,F\subseteq\N\),
\[
 \left\{b\in M^y:
 \varphi^M(a_n,b)\leq r\ (n\in E),\quad
 \varphi^M(a_n,b)\geq s\ (n\in F)\right\}
 \neq\varnothing.
\]
The formula is \emph{NIP on \(M\)} if it lacks IP on \(M\).

The formula \(\varphi\) is \emph{stable on \(M\)} if, for all sequences
\((a_m)\) in \(M^x\) and \((b_n)\) in \(M^y\) for which both iterated
limits exist, one has
\[
 \lim_m\lim_n\varphi^M(a_m,b_n)
 =
 \lim_n\lim_m\varphi^M(a_m,b_n).
\]
This is Grothendieck's double-limit condition
\cite[Definition~4.7]{KhanakiBanach}. Stability on \(M\) implies NIP on \(M\)
\cite[Fact~4.10]{KhanakiBanach}.

Let \(T\) be a complete theory. A fixed formula \(\varphi\) is
\emph{stable in \(T\)} if it is stable on every model of \(T\). It has
\emph{IP in \(T\)} if it has IP on some model of \(T\), and is
\emph{NIP in \(T\)} if it is NIP on every model of \(T\).
The theory \(T\) is NIP if every formula is NIP in \(T\).

We call the norm formula \(\lVert x+y\rVert\) stable or NIP on a Banach
space \(X\) when the corresponding property holds on every pair of ball sorts.
For stability, this is the Banach-space double-limit condition on bounded
sequences \cite{KrivineMaurey}.
It has IP on \(X\) if it has IP on some pair of ball sorts. By homogeneity,
the two unit-ball sorts suffice to test these properties. Given a common
positive bound \(R\) on the norms of the vectors in a witness, dividing
the vectors by \(R\) scales all norm values, and any IP thresholds, by
\(1/R\).

Khanaki proved that the norm formula is NIP on the Figiel--Johnson version of
Tsirelson space
\cite[Corollary~4.19]{KhanakiBanach}. The same space is not stable
\cite[Corollary~4.14]{KhanakiBanach}, by the Krivine--Maurey theorem and the
fact that it contains no copy of any
\(\ell_p\) \cite[Section~2]{FJ74}.
Khanaki and Pillay discuss the question whether the same formula is NIP
in the complete theory of that space
\cite[Introduction]{KhanakiPillay}.

Let \(Y\) be \(H\) with the norm \(\sqrt{h^2+q^2}\).
By Theorem~\ref{thm:y}, the norm formula is unstable on \(Y\)
and NIP in \(\Th(Y)\). Lupini uses a dense linear order to index the
coordinates of an additive Hilbert renorming \cite[Corollary~7.4]{Lupini2026Norm}
and proves that its complete theory is NIP \cite[Theorem~9.5]{Lupini2026Norm}.
His seminorm includes the individual coordinate moduli as separate terms
(Section~\ref{subsec:y-definition}). In the supremum defining \(q\) for a
finitely supported pair \((u,v)\), each \(u_i\) can be paired with
a zero \(v_j\) by taking \(j>i\) beyond the support of \(v\). A coordinate
\(v_j\) can be paired with a zero \(u_i\) only if one occurs before \(j\).
For the first unit vector \(e_1\) of \(\ell_2(\N)\), we have
\(q(0,e_1)=0\) because no index precedes \(1\).
Lemma~\ref{lem:sparse-trace} includes a finite initial segment to determine
which second-copy coordinates admit such a pair.

We say that a Banach space is \(C\)-isomorphic to Hilbert space when it admits
an isomorphism \(J\) onto a Hilbert space with
\(\lVert J\rVert\lVert J^{-1}\rVert\leq C\).
Section~\ref{subsec:x-definition} defines \(X\) as the \(\ell_2\)-sum
of finite-dimensional blocks whose quadratic norms depend on membership
in finite sets. Lupini uses the same membership relation with additive
block norms \cite[Proposition~8.3]{Lupini2026Norm}.

\begin{theorem}\label{thm:main}
There are separable Banach spaces \(X\) and \(Y\), each
\(\sqrt3\)-isomorphic to Hilbert space, with the following properties.
\begin{enumerate}[label=\textup{(\roman*)},leftmargin=2.2em]
\item The norm formula is stable on \(X\) and has IP in \(\Th(X)\).
\item The norm formula is unstable on \(Y\) and is NIP in \(\Th(Y)\).
\end{enumerate}
Both spaces are uniformly convex, and their norms satisfy
\[
 \lVert x\rVert,\lVert y\rVert\leq1
 \quad\Longrightarrow\quad
 \left\lVert\frac{x+y}{2}\right\rVert^2
 \leq1-\frac1{12}\lVert x-y\rVert^2.
\]
There is also a separable elementary extension \(W\) of \(X\) satisfying the
same estimate and such that the norm formula has IP on \(W\). The space \(W\) is
\(\sqrt3\)-isomorphic to Hilbert space and contains no isomorphic copy of
\(c_0\).
\end{theorem}

Table~\ref{tab:combinations} gives the stability and independence
properties of the norm formula for \(\ell_2\), \(X\), \(Y\),
\(X\oplus_2Y\), and \(W\). For the sum, these are proved in
Corollary~\ref{cor:xy}. All five spaces are separable and
\(\sqrt3\)-isomorphic to \(\ell_2\).

\begin{table}[htbp]
\centering
\renewcommand{\arraystretch}{1.2}
\begin{tabular}{c@{\qquad}cc}
\toprule
 & NIP in the complete theory & IP in the complete theory\\
\midrule
stable on the space & \(\ell_2\) & \(X\)\\
unstable and NIP on the space & \(Y\) & \(X\oplus_2Y\)\\
IP on the space & --- & \(W\)\\
\bottomrule
\end{tabular}
\caption{The norm formula has the indicated properties on each space and in its complete theory.}
\label{tab:combinations}
\end{table}

IP on a space implies IP in its complete theory, so the last row can
occur only in the IP column.

The norm formula is stable both on Hilbert space with its usual norm and
in its complete theory. For Lupini's additive block norm, the formula is
stable on the space and has IP in its complete theory
\cite[Proposition~8.3]{Lupini2026Norm}.
The norm formula of his dense-order renorming is NIP in its complete
theory \cite[Corollary~7.4]{Lupini2026Norm}.

By Lemma~\ref{lem:local-nip-sum}, the norm formula is NIP on the
\(\ell_2\)-sum of Lupini's block and dense-order spaces. The coordinate
vectors in the dense-order summand give distinct iterated limits
using the seminorm values in \cite[eq.~(5)]{Lupini2026Norm},
so the formula is unstable on the sum. The block summand has arbitrarily
large finite independence patterns with fixed thresholds
\cite[Proposition~8.3]{Lupini2026Norm}. The same patterns occur in the sum,
so Lemma~\ref{lem:finite-pattern-transfer}(i) gives IP in its complete theory.
An elementary extension of Lupini's stable example that realizes an
IP witness gives a space on which the norm formula has IP.
By downward L\"owenheim--Skolem
\cite[Proposition~7.3]{BBHU08}, this extension can be chosen separable,
as in the proof of Theorem~\ref{thm:hilbertian-ip-space}.

These examples realize the same five combinations as the spaces in
Table~\ref{tab:combinations}. The five spaces in the table also satisfy
the common midpoint estimate in Theorem~\ref{thm:main}.

A Banach space is \emph{superstable} if every Banach
space finitely representable in it is stable \cite{Raynaud}. Equivalently,
all its ultrapowers are stable by \cite[Theorem~6.3]{Heinrich}, since
stability passes to isometric subspaces. The finite membership
patterns constructed in Section~\ref{subsec:x-definition} give IP on
every ultrapower of \(X\) by a nonprincipal ultrafilter on \(\N\), by
Lemma~\ref{lem:finite-pattern-transfer}(i). The same witnesses give IP
on the corresponding ultrapowers of \(X\oplus_2Y\) in the proof of
Corollary~\ref{cor:xy}. Since IP implies instability, both spaces fail
superstability.

\section{Preliminaries}

We use bounded real-valued formulas on fixed sorts. To apply a result stated
for \([0,1]\)-valued formulas to a bounded real-valued formula, we compose the
formula with a positive affine map sending its range into \([0,1]\). This
preserves stability, IP, and NIP.

\subsection{Stability and pointwise convergence on local type spaces}

For a bounded formula \(\varphi(x,y)\) on a metric structure \(M\), put
\(\widetilde\varphi(y,x)=\varphi(x,y)\), and let
\(S_{\widetilde\varphi}(M)\) be the compact space of complete
\(\widetilde\varphi\)-types over \(M\)
\cite[Definition~6.1 and Fact~6.2]{BU10Author}. If \(b\in N^y\) realizes a type
\(p\) in an elementary extension \(N\succeq M\), then \(p\) is given by the
conditions \(\varphi(a,y)=\varphi^N(a,b)\) for all \(a\in M^x\).
For each \(a\in M^x\), the evaluation function
\[
 \varphi(a,\cdot):S_{\widetilde\varphi}(M)\longrightarrow\R
\]
sends \(p\) to \(\varphi^N(a,b)\) for any realization \(b\) of \(p\) in
\(N\succeq M\). The topology is pointwise convergence of the functions
\(a\mapsto\varphi(a,p)\) on \(M^x\). By Grothendieck's double-limit
criterion \cite[Theorem~6]{Gro52}, stability on \(M\) is equivalent to
relative weak compactness of the family of evaluation functions in
\(C(S_{\widetilde\varphi}(M))\)
\cite[Definition~4.7 and Remark~4.8(i)]{KhanakiBanach}.
The Eberlein--\v Smulian theorem then gives the implication from stability
to NIP \cite[Fact~4.10]{KhanakiBanach}.

When \(M\) is separable and the language is countable,
\(S_{\widetilde\varphi}(M)\) is compact and metrizable. On this space,
pointwise subsequence convergence characterizes NIP on \(M\) by the
Bourgain--Fremlin--Talagrand theorem \cite{BFT}, in the form stated by
Khanaki \cite[Lemma~4.2(v)]{KhanakiBanach}.

\begin{lemma}\label{lem:sequential-nip}
Let \(M\) be a separable metric structure in a countable language, and let
\(\varphi(x,y)\) be a bounded formula on fixed sorts. Then \(\varphi\) is
NIP on \(M\) if and only if every sequence \(\varphi(a_n,\cdot)\), with
\(a_n\in M^x\), has a pointwise convergent subsequence on
\(S_{\widetilde\varphi}(M)\).
\end{lemma}

The convergence is pointwise on the type space. Equivalently, if
\(M^*\succeq M\) realizes every
\(\widetilde\varphi\)-type over \(M\), as happens when \(M^*\) is
\(\lvert M\rvert^+\)-saturated, the subsequence may be chosen so that
\(\varphi^{M^*}(a_n,b)\) converges for every \(b\in(M^*)^y\).

\subsection{Finite boxes}

For \(n\geq1\), write
\([n]=\{1,\ldots,n\}\).

\begin{definition}\label{def:finite-ip-pattern}
Let \(\psi(x,y)\) be a bounded formula on fixed sorts of a structure \(M\),
and let \(\bar a=(a_1,\ldots,a_n)\in(M^x)^n\). Its \emph{evaluation set} is
\[
 Q_{\psi,\bar a}^M
 =\left\{\bigl(\psi^M(a_1,b),\ldots,\psi^M(a_n,b)\bigr):b\in M^y\right\}
 \subseteq\R^n.
\]
For \(r<s\), we say that this evaluation set \emph{determines an
\(n\)-dimensional box between \(r\) and \(s\)} if, for every
\(S\subseteq[n]\), there is a parameter \(b_S\in M^y\) such that
\[
 \psi^M(a_i,b_S)\geq s\quad(i\in S),
 \qquad
 \psi^M(a_i,b_S)\leq r\quad(i\notin S).
\]
\end{definition}

This is the form of Ben Yaacov's box definition
\cite[Definition~2.3]{BenYaacovVC} with common thresholds in every
coordinate, extended to bounded real-valued formulas. The box has
width \(s-r\).

\begin{lemma}\label{lem:finite-pattern-transfer}
Let \(\psi(x,y)\) be a bounded formula on fixed sorts of a structure \(M\).
\begin{enumerate}[label=\textup{(\roman*)},leftmargin=2.2em]
\item If, for every \(n\), some evaluation set in \(M\) determines an
\(n\)-dimensional box between fixed thresholds \(r<s\), then \(\psi\) has
IP on \(M^{\mathcal U}\) for every nonprincipal ultrafilter
\(\mathcal U\) on \(\N\).
\item If \(\psi\) has IP in \(\Th(M)\), witnessed with thresholds \(r<s\),
then, for every \(n\) and every \(0<\delta<(s-r)/2\), some evaluation set in
\(M\) determines an \(n\)-dimensional box between \(r+\delta\) and
\(s-\delta\).
\end{enumerate}
\end{lemma}

\begin{proof}
For (i), let \(\mathcal U\) be a nonprincipal ultrafilter on \(\N\).
For each \(n\), choose tuples \(a_i^{(n)}\) and parameters \(b_S^{(n)}\)
witnessing the \(n\)-dimensional box. Represent \(a_i\) by \(a_i^{(n)}\)
for \(n\geq i\), using a fixed element at the remaining coordinates.
For finite disjoint \(E,F\subseteq\N\), represent the parameter by
\(b_{F\cap[n]}^{(n)}\). The required inequalities hold for all sufficiently
large \(n\), and hence in the ultrapower by the ultraproduct theorem
\cite[Theorem~5.4]{BBHU08}.

For (ii), write \(t_+=\max\{t,0\}\). The sentence
\[
\inf_{\substack{u_1,\ldots,u_n\\(v_S)_{S\subseteq[n]}}}
\max_{S\subseteq[n]}
\max\left\{
 \max_{i\in S}\bigl(s-\psi(u_i,v_S)\bigr)_+,
 \max_{i\notin S}\bigl(\psi(u_i,v_S)-r\bigr)_+
\right\}
\]
takes the infimum of the largest violation of the box inequalities.
A maximum over the empty set is \(0\), and the quantifiers range over
the fixed sorts of \(\psi\). The sentence has value zero in
a model witnessing IP and, by completeness, in \(M\). Choose tuples in
\(M\) for which the maximum in the display is less than \(\delta\).
Then \(\psi(u_i,v_S)\) is at least \(s-\delta\) for \(i\in S\) and at
most \(r+\delta\) for \(i\notin S\).
\end{proof}

A continuous connective is a continuous map from a finite product of
compact real intervals to a compact real interval, applied pointwise to
formula values. NIP in a theory is preserved by continuous connectives
\cite[Proposition~3.7]{BenYaacovVC}.

\begin{lemma}\label{lem:nip-connectives}
Let \(T\) be a complete continuous theory. For \(1\leq j\leq d\), let
\(I_j\) be a compact real interval and let \(\psi_j(x,y)\) be a formula
with values in \(I_j\) that is NIP in \(T\). Let \(J\) be a compact real
interval and let \(F:I_1\times\cdots\times I_d\to J\) be continuous. Then
\[
                 F(\psi_1(x,y),\ldots,\psi_d(x,y))
\]
is NIP in \(T\).
\end{lemma}

\begin{proof}
A bounded formula is NIP in \(T\) if and only if, for every indiscernible
sequence \((a_i)\) and every parameter \(b\), its values at \((a_i,b)\) converge
\cite[Lemma~5.4]{BenYaacovVC}; see also
\cite[Definition~4.1]{Anderson}. Fix such a sequence and parameter. Each
sequence \(\psi_j(a_i,b)\) converges. Continuity of \(F\) therefore gives
convergence of \(F(\psi_1(a_i,b),\ldots,\psi_d(a_i,b))\), and the same
criterion gives NIP.
\end{proof}

\subsection{Hilbert norm comparisons and uniform convexity}

The quadratic combination of a Hilbert norm and a seminorm bounded by
\(\sqrt2\) times the Hilbert norm satisfies a midpoint estimate for
uniform convexity.

\begin{lemma}\label{lem:uniform-convexity}
Let \(K\) be a Hilbert space and let \(p\) be a seminorm satisfying
\(p(z)\leq\sqrt2\lVert z\rVert_2\). Define
\[
                    N(z)=\sqrt{\lVert z\rVert_2^2+p(z)^2}.
\]
Then
\[
\lVert z\rVert_2\leq N(z)\leq\sqrt3\lVert z\rVert_2,
\]
and, if \(N(x),N(y)\leq1\), then
\[
N\!\left(\frac{x+y}{2}\right)^2
\leq1-\frac1{12}N(x-y)^2.
\]
In particular, \(N\) is uniformly convex, with a modulus independent of
the dimension of \(K\) and of the seminorm \(p\).
\end{lemma}

\begin{proof}
Subadditivity of \(\lVert\cdot\rVert_2\) and \(p\), followed by Minkowski's
inequality in \(\R^2\), shows that \(N\) is a norm. The norm equivalence
follows directly from \(p(z)\leq\sqrt2\lVert z\rVert_2\). For the convexity
estimate, convexity of
\(p\), the scalar
inequality \(((a+b)/2)^2\leq(a^2+b^2)/2\), and the Hilbert parallelogram
identity give
\begin{align*}
N\!\left(\frac{x+y}{2}\right)^2
&\leq
\frac{\lVert x\rVert_2^2+\lVert y\rVert_2^2}{2}
-\frac14\lVert x-y\rVert_2^2
+\frac{p(x)^2+p(y)^2}{2}\\
&\leq1-\frac14\lVert x-y\rVert_2^2
\leq1-\frac1{12}N(x-y)^2. \qedhere
\end{align*}
\end{proof}

\section{Stable examples}

\subsection{Hilbert space}\label{subsec:hilbert-stability}

The complete theory of infinite-dimensional
Hilbert spaces is stable \cite[Example~6.5(iii)]{KhanakiBanach}.
Consequently, the norm formula is stable, and hence NIP, in
\(\Th(\ell_2)\).
The norm formula is also stable on every finite-dimensional Hilbert space,
by restriction from \(\ell_2\).

\subsection{Stability of \texorpdfstring{\(\ell_2\)}{l2}-sums of
finite-dimensional spaces}

Let
\[
       Z=\left(\bigoplus_{k=1}^{\infty}E_k\right)_{\ell_2},
\]
where each \(E_k\) is a finite-dimensional Banach space. Thus
\[
       \lVert z\rVert_Z^2=\sum_k\lVert z_k\rVert_{E_k}^2.
\]

We use the elementary inequality
\begin{equation}\label{eq:defect}
 \left|
 \lVert u+v\rVert^2-\lVert u\rVert^2-\lVert v\rVert^2
 \right|
 \leq 2\lVert u\rVert\,\lVert v\rVert,
\end{equation}
valid in every normed space.

\begin{lemma}\label{lem:pythagoras}
Suppose that \((z_m)\) is a bounded weakly null sequence in \(Z\) and
\(\lVert z_m\rVert\to\lambda\).  Then for every \(u\in Z\),
\[
       \lVert u+z_m\rVert^2
       \longrightarrow
       \lVert u\rVert^2+\lambda^2.
\]
\end{lemma}

\begin{proof}
For each \(m,k\), put
\[
D_{m,k}=\lVert u_k+z_{m,k}\rVert^2-\lVert u_k\rVert^2-\lVert z_{m,k}\rVert^2.
\]
For each fixed \(k\), weak convergence gives \(z_{m,k}\to0\) in norm because
\(E_k\) is finite-dimensional. Thus \(D_{m,k}\to0\), and the sum of
\(D_{m,k}\) over \(k\leq N\) tends to zero for each fixed \(N\). For the tail,
\eqref{eq:defect} and the Cauchy--Schwarz inequality give
\begin{align*}
\sum_{k>N}|D_{m,k}|
&\leq 2\sum_{k>N}\lVert u_k\rVert\lVert z_{m,k}\rVert\\
&\leq
2\left(\sum_{k>N}\lVert u_k\rVert^2\right)^{1/2}
\sup_m\lVert z_m\rVert.
\end{align*}
This bound tends to zero as \(N\to\infty\), uniformly in \(m\). Hence
\(\lVert u+z_m\rVert^2-\lVert u\rVert^2-\lVert z_m\rVert^2
=\sum_k D_{m,k}\to0\). Since \(\lVert z_m\rVert^2\to\lambda^2\), the
asserted limit follows.
\end{proof}

\begin{theorem}\label{thm:local-stability}
For every sequence of finite-dimensional Banach spaces \((E_k)\), the norm
formula is stable on
\[
             \left(\bigoplus_kE_k\right)_{\ell_2}.
\]
\end{theorem}

\begin{proof}
Let \((x_m)\) and \((y_n)\) be bounded sequences and suppose the two iterated
limits of \(\lVert x_m+y_n\rVert\) exist. The space is reflexive, since it is
an \(\ell_2\)-sum of finite-dimensional spaces. Pass to weakly convergent
subsequences and refine them so that
\[
x_m\rightharpoonup x,\quad \lVert x_m-x\rVert\to\lambda,
\qquad
y_n\rightharpoonup y,\quad \lVert y_n-y\rVert\to\mu.
\]
Each inner limit is unchanged on the chosen subsequence, and the outer
limit is then the limit of a subsequence of the original scalar sequence.
For each fixed \(m\), apply Lemma~\ref{lem:pythagoras} with \(u=x_m+y\)
and the weakly null sequence \(y_n-y\). This gives
\[
\lim_n\lVert x_m+y_n\rVert^2=\lVert x_m+y\rVert^2+\mu^2.
\]
Apply the lemma again with \(u=x+y\) and \(x_m-x\) to get
\(\lim_m\lVert x_m+y\rVert^2=\lVert x+y\rVert^2+\lambda^2\).
Taking the outer limit and then the square root gives
\[
\lim_m\lim_n\lVert x_m+y_n\rVert
=\left(\lVert x+y\rVert^2+\lambda^2+\mu^2\right)^{1/2}.
\]
Interchanging the two sequences gives the same value for the reverse
iterated limit. This proves stability by the double-limit condition.
\end{proof}

The theorem is a special case of the classical fact that, for
\(1\leq p<\infty\), an \(\ell_p\)-sum of stable Banach spaces is stable; see
Baudier
\cite[discussion following Example~3.10]{Baudier}.
Lupini proves stability of \(\ell_2\)-sums of finite-dimensional spaces
using ultrafilter limits \cite[Lemma~8.2]{Lupini2026Norm}.
Lemma~\ref{lem:pythagoras} gives the convergent-norm form of property
\((m_2)\) of Kalton and Werner \cite[Definition~3.1(a)]{KaltonWerner}.
A space has this property if every weakly null sequence \((z_m)\) and
every vector \(u\) satisfy
\[
\limsup_m\lVert u+z_m\rVert^2
=\lVert u\rVert^2+\limsup_m\lVert z_m\rVert^2.
\]

\subsection{An \texorpdfstring{\(\ell_2\)}{l2}-sum of \texorpdfstring{\(\ell_\infty^n\)}{linfinity n}}

Consider the elementary example
\[
       X_\infty=\left(\bigoplus_{n=1}^{\infty}\ell_\infty^n\right)_{\ell_2}.
\]
Theorem~\ref{thm:local-stability} shows that the norm formula is stable on
\(X_\infty\). In the \(n\)-th summand put
\[
       a_i=e_i\quad(1\leq i\leq n),
       \qquad b_S=\mathbf1_S=\sum_{j\in S}e_j.
\]
Then
\[
 \lVert a_i+b_S\rVert=
 \begin{cases}
 2,&i\in S,\\
 1,&i\notin S.
 \end{cases}
\]
Thus the evaluation sets determine boxes between \(1\) and \(2\) in every
dimension, and Lemma~\ref{lem:finite-pattern-transfer}(i) gives IP in
\(\Th(X_\infty)\).

The corresponding ultrapowers also contain \(c_0\) isometrically. Let
\(\mathcal U\) be a nonprincipal ultrafilter on \(\N\), put \(x_i^{(n)}=e_i\)
in the \(n\)-th summand for \(n\geq i\), and put \(x_i^{(n)}=0\) for \(n<i\).
For scalars \(\alpha_1,\ldots,\alpha_j\), the norm of
\(\sum_{i=1}^j\alpha_i x_i^{(n)}\) is \(\max_{1\leq i\leq j}|\alpha_i|\)
for all \(n\geq j\). Thus
\[
\left\lVert\sum_{i=1}^j\alpha_i[x_i^{(n)}]_{\mathcal U}\right\rVert
=\max_{1\leq i\leq j}|\alpha_i|.
\]
The closed span of these ultrapower vectors is isometric to \(c_0\).

\subsection{Stability on \texorpdfstring{\(X\)}{X} and IP in
\texorpdfstring{\(\Th(X)\)}{Th(X)}}\label{subsec:x-definition}

For each \(n\geq1\), let \(H_n\) be the Hilbert space with
orthonormal basis
\[
       f_1,\ldots,f_n,
       \qquad (e_S)_{S\subseteq[n]}.
\]
For \(z\in H_n\), put
\[
p_n(z)=\max_{\substack{S\subseteq[n]\\ i\in S}}
       \bigl(|\langle z,f_i\rangle|+|\langle z,e_S\rangle|\bigr),
\qquad
N_n(z)=\left(\lVert z\rVert_2^2+p_n(z)^2\right)^{1/2}.
\]
Since \(p_n(z)\leq\sqrt2\lVert z\rVert_2\),
\begin{equation}\label{eq:x-equivalence}
              \lVert z\rVert_2\leq N_n(z)\leq\sqrt3\lVert z\rVert_2
\end{equation}
uniformly in \(n\). If \(i\in S\), the pair \((i,S)\) gives
\(p_n(f_i+e_S)=2\). If \(i\notin S\), every sum in the definition is at
most \(1\), and the pair \((i,\{i\})\) attains \(1\). This includes
\(S=\varnothing\). Therefore
\begin{equation}\label{eq:x-cubes}
 N_n(f_i+e_S)=
 \begin{cases}
 \sqrt6,&i\in S,\\
 \sqrt3,&i\notin S.
 \end{cases}
\end{equation}
The vectors \(f_i\) and \(e_S\) have \(N_n\)-norm at most
\(\sqrt2\).

Lupini uses the same membership relation in finite-dimensional blocks
\cite[Proposition~8.3]{Lupini2026Norm}. His seminorm is the maximum
of the absolute values of individual coordinates and the absolute
values of coordinate sums at pairs \((i,S)\) with \(i\in S\).
He adds a multiple of this seminorm to the Hilbert norm of each block.

Define
\[
              X=\left(\bigoplus_{n=1}^{\infty}(H_n,N_n)\right)_{\ell_2}.
\]

\begin{theorem}\label{thm:x}
The norm formula is stable on \(X\) and has IP in
\(\Th(X)\).
\end{theorem}

\begin{proof}
Stability on \(X\) follows from Theorem~\ref{thm:local-stability}. In the
\(n\)-th summand, the vectors and parameters above lie in fixed balls, and
\eqref{eq:x-cubes} shows that the corresponding evaluation set determines
a box between \(\sqrt3\) and \(\sqrt6\).
Lemma~\ref{lem:finite-pattern-transfer}(i) therefore gives IP in the complete
theory.
\end{proof}

On the underlying Hilbert sum define
\[
 p_{\oplus}(z)=\left(\sum_n p_n(z_n)^2\right)^{1/2}.
\]
The bound on \(p_n\) makes this expression finite. Coordinatewise
subadditivity and Minkowski's inequality in \(\ell_2\) show that
\(p_{\oplus}\) is a seminorm. The uniform equivalence
\eqref{eq:x-equivalence} identifies the underlying vector space of \(X\)
with the Hilbert direct sum, and summing the block identities gives
\[
 \lVert z\rVert_X
 =\left(\lVert z\rVert_2^2+p_{\oplus}(z)^2\right)^{1/2},
 \qquad
 p_{\oplus}(z)\leq\sqrt2\lVert z\rVert_2.
\]
Lemma~\ref{lem:uniform-convexity} therefore gives the equivalence and the
uniform convexity estimate asserted for \(X\) in
Theorem~\ref{thm:main}.

Both \(X_\infty\) and \(X\) are reflexive and stable. The copy of \(c_0\)
in \(X_\infty^{\mathcal U}\) rules out superreflexivity of \(X_\infty\)
\cite[Proposition~6.4]{Heinrich}. The uniform convexity of \(X\) gives
superreflexivity. The boxes in the two block sums give IP on every
ultrapower of either space by a nonprincipal ultrafilter on \(\N\),
by Lemma~\ref{lem:finite-pattern-transfer}(i). Neither space is superstable.

\section{Renormings defined by \texorpdfstring{\(i<j\)}{i<j}}

\subsection{Bollob\'as's norm and the quadratic renorming}
\label{subsec:y-definition}

Let
\[
H=\ell_2(\N)\oplus_2\ell_2(\N),\qquad h(z)=\lVert z\rVert_2
\]
and define
\begin{equation}\label{eq:q-definition}
q(u,v)=\sup_{i<j}|u_i+v_j|,
\qquad
\lVert(u,v)\rVert_Y
=\left(\lVert u\rVert_2^2+\lVert v\rVert_2^2+q(u,v)^2\right)^{1/2}.
\end{equation}
We write \(q(z)=q(u,v)\) for \(z=(u,v)\). An index pair \((i,j)\) is
\emph{admissible} if \(i<j\). The supremum of the absolute values of
the corresponding linear functionals defines a seminorm, and
\[
q(u,v)\leq\sqrt2\lVert(u,v)\rVert_2,
\]
and therefore
\begin{equation}\label{eq:y-equivalence}
\lVert z\rVert_2\leq\lVert z\rVert_Y
\leq\sqrt3\lVert z\rVert_2.
\end{equation}
Let \(Y=(H,\lVert\cdot\rVert_Y)\). For \(x_i=(e_i,0)\) and
\(y_j=(0,e_j)\), one has
\[
q(x_i+y_j)=
\begin{cases}
2,&i<j,\\
1,&i\geq j.
\end{cases}
\]
If \(i<j\), the admissible pair \((i,j)\) contains both nonzero
coordinates. If \(i\geq j\), each admissible pair contains at most one
nonzero coordinate, and \((i,i+1)\) contains the first and gives the value \(1\).
Hence
\begin{equation}\label{eq:half-graph}
\lVert x_i+y_j\rVert_Y=
\begin{cases}
\sqrt6,&i<j,\\
\sqrt3,&i\geq j.
\end{cases}
\end{equation}
Thus \(\lim_i\lim_j\lVert x_i+y_j\rVert_Y=\sqrt6\), while
\(\lim_j\lim_i\lVert x_i+y_j\rVert_Y=\sqrt3\). The norm formula is therefore
unstable on \(Y\).

Applying Lemma~\ref{lem:uniform-convexity} to the Hilbert space \(H\) and
the seminorm \(q\) gives the uniform convexity estimate for \(Y\) stated
in Theorem~\ref{thm:main}.

The seminorm \(q\) also occurs in Bollob\'as's unstable renorming of
\(\ell_2\), described by Garling \cite{Garling1982} and reproduced by Iovino
\cite[remark following Theorem~8.7 in the author version]{IovinoStableBanachI}.
Iovino uses the norm
\[
 \lVert z\rVert_\eps
 =\lVert z\rVert_2+\eps\sup_{i<j}|z_{2i}+z_{2j-1}|,
 \qquad \eps>0.
\]
Under the identification \(u_i=z_{2i}\) and \(v_j=z_{2j-1}\), this norm
is \(h+\eps q\).

Let \(I=\mathbb Q\cap(0,1)\) with its usual order. On
\(H_I=\ell_2(I)\oplus_2\ell_2(I)\) Lupini defines
\[
 \gamma(a,b)=\max\left\{
     \sup_{d\in I}|a_d|,\ \sup_{d\in I}|b_d|,\ \sup_{d<e}|a_d+b_e|
     \right\}
\]
and adds a positive multiple of \(\gamma\) to the Hilbert norm
\cite[Corollary~7.4]{Lupini2026Norm}. The resulting renorming has NIP
complete theory by a proof in an expansion with an ordered orthonormal
family \cite[Theorem~9.5]{Lupini2026Norm}.
In \(q\), no admissible pair has second index \(1\), so \(q(0,e_1)=0\).
In \(\gamma\), the second supremum includes every coordinate,
giving \(\gamma(0,e_d)=1\) for each \(d\in I\).

\subsection{Sparse approximation and finite fat-shattering dimension}
\label{subsec:fat-shattering}

\begin{definition}\label{def:fat-shattering}
Let \(D\) be a set, let \(\mathcal F\subseteq\R^D\), and let
\(\eps>0\). Points \(x_1,\ldots,x_n\in D\) are
\(\eps\)-fat-shattered by \(\mathcal F\) if there are real numbers
\(c_1,\ldots,c_n\) such that, for every \(S\subseteq[n]\), some
\(f\in\mathcal F\) satisfies
\[
f(x_i)\geq c_i+\eps\quad(i\in S),
\qquad
f(x_i)\leq c_i-\eps\quad(i\notin S).
\]
\end{definition}

The \(\eps\)-fat-shattering dimension of \(\mathcal F\) is the supremum of
the integers \(n\) for which such a family of \(n\) points exists. It is
infinite if these integers are unbounded
\cite[Definition~2.2]{AlonBenDavidCesaBianchiHaussler1997}.

Let \(\psi(x,y)\) be a bounded formula on fixed sorts of \(M\), and let
\(\mathcal F=\{\psi^M(\cdot,b):b\in M^y\}\) be the corresponding class of
functions on \(M^x\). If the evaluation set for \(a_1,\ldots,a_n\)
determines a box between \(r<s\) in the sense of
Definition~\ref{def:finite-ip-pattern}, then \(\mathcal F\) fat-shatters
these points at scale \((s-r)/2\), with each threshold equal to \((r+s)/2\).

If \(\psi\) has IP in \(\Th(M)\) with thresholds \(r<s\), fix
\(\delta=(s-r)/4\). Lemma~\ref{lem:finite-pattern-transfer}(ii) gives boxes
in \(M\) between \(r+\delta\) and \(s-\delta\) in every dimension, so
\(\mathcal F\) has infinite fat-shattering dimension at scale
\((s-r)/4\). Thus if \(\mathcal F\) has finite fat-shattering dimension at
every positive scale, then \(\psi\) is NIP in \(\Th(M)\).

Lupini bounds the fat-shattering dimension of the translates of \(\gamma\)
by truncating small coordinates, counting the possible orders of the
retained indices, and counting sign patterns of polynomials
\cite[Lemma~7.2]{Lupini2026Norm}. For \(q\), we also truncate small
coordinates. For fixed, finitely supported test vectors \(x_1,\ldots,x_n\)
and a finitely supported parameter \(y\), Lemma~\ref{lem:sparse-trace}
expresses each \(q(x_r+y)\) as a finite maximum that includes admissible
pairs with a zero coordinate. To count threshold patterns as \(y\) varies,
fix the copy of each nonzero coordinate of \(y\) and the positions of the
parameter indices relative to one another and to the indices in the
test-vector supports and the initial segment used in the finite maximum.
Write \(t_i\) for the value of the \(i\)th nonzero parameter coordinate.
Each nonconstant comparison with a threshold tests \(t_i\) or
\(t_i+t_j\) against cut values determined by the thresholds and the
test-vector coordinates. The cut values partition \(\R\) into singletons
and open intervals. The threshold pattern is determined by the piece
containing each \(t_i\) and each \(t_i+t_j\).

Let \(K\) be a Hilbert space with a fixed orthonormal basis, and write
\(\operatorname{supp}z\) for the support of \(z\) in this basis.  Let
\(p:K\to[0,\infty)\), and suppose that for every \(\eta>0\) there is a map
\(P_\eta:K\to K\), a nonnegative integer-valued function \(a(C,\eta)\)
defined for \(C\geq0\), and a number
\(\omega(\eta)\), with \(\omega(\eta)\to0\) as \(\eta\downarrow0\), such that
\begin{align}
 |\operatorname{supp}(P_\eta z)|
 &\leq a(C,\eta) &&(\lVert z\rVert_2\leq C),
 \label{eq:criterion-support}\\
 |p(x+y)-p(P_\eta x+P_\eta y)|
 &\leq\omega(\eta) &&(x,y\in K).
 \label{eq:criterion-approximation}
\end{align}
Assume also that, for every pair of nonnegative integers \(m,k\), there are
constants \(A_{m,k}\geq1\) and \(d_{m,k}\geq0\) with the following
property.  If \(x_1,\ldots,x_n\) have support size at most \(m\) and
\(c_1,\ldots,c_n\in\R\), define the \emph{threshold trace} of a parameter
\(y\) by
\begin{equation}\label{eq:criterion-traces}
 \left(\mathbf1_{\{p(x_i+y)\geq c_i\}}\right)_{i=1}^n.
\end{equation}
As \(y\) ranges over vectors of support size at most \(k\), the number of
distinct threshold traces is at most
\begin{equation}\label{eq:criterion-trace-bound}
                         A_{m,k}(n+1)^{d_{m,k}}.
\end{equation}

\begin{proposition}
\label{prop:sparse-norming-criterion}
Suppose that \(p\) satisfies the approximation hypotheses
\eqref{eq:criterion-support}--\eqref{eq:criterion-approximation} and the
polynomial trace bound \eqref{eq:criterion-trace-bound}. Then, for every
\(0\leq R_x,R_y<\infty\) and \(\eps>0\), the class
\[
 \mathcal P_{R_x,R_y}
 =\{f_y:R_x\Ball_K\to\R:f_y(x)=p(x+y),\ \lVert y\rVert_2\leq R_y\}
\]
has finite \(\eps\)-fat-shattering dimension.
\end{proposition}

\begin{proof}
Choose \(\eta>0\) so that \(\omega(\eta)<\eps\), and set
\[
 m=a(R_x,\eta),\qquad k=a(R_y,\eta).
\]
Suppose \(x_1,\ldots,x_n\) are \(\eps\)-fat-shattered, with thresholds
\(c_1,\ldots,c_n\), and choose a parameter \(y_S\) for every
\(S\subseteq[n]\).  By \eqref{eq:criterion-support}, the vectors
\(P_\eta x_i\) and \(P_\eta y_S\) have support at most \(m\) and \(k\),
respectively.  Since the error in \eqref{eq:criterion-approximation} is
strictly smaller than the margin, \(P_\eta y_S\) has threshold trace
\(\mathbf1_S\) on the vectors \(P_\eta x_i\). Hence
\begin{equation}\label{eq:criterion-exponential-polynomial}
                         2^n\leq A_{m,k}(n+1)^{d_{m,k}}.
\end{equation}

The left-hand side grows exponentially in \(n\), while the right-hand side
is bounded by a polynomial in \(n\). Hence the inequality fails for all
sufficiently large \(n\), with the bound depending only on \(R_x,R_y,\eps\).
This proves finiteness of the fat-shattering dimension.
\end{proof}

For \(0\leq R_x,R_y<\infty\), set
\[
\mathcal Q_{R_x,R_y}
=\{f_y:R_x\Ball_H\to\R:f_y(x)=q(x+y),\ \lVert y\rVert_2\leq R_y\}.
\]
Index coordinate \(j\) in copy \(\iota\) by \((\iota,j)\), so that the
support of \(z=(u,v)\in H\) is
\[
 \operatorname{supp}z
 =\bigl(\{1\}\times\operatorname{supp}u\bigr)
  \cup\bigl(\{2\}\times\operatorname{supp}v\bigr)
 \subseteq\{1,2\}\times\N.
\]
For \(\eta>0\), let \(T_\eta z\) retain the coordinates of \(z\) having
absolute value at least \(\eta\). If \(C\geq0\) and
\(\lVert z\rVert_2\leq C\), then
\begin{equation}\label{eq:support-bound}
 |\operatorname{supp}(T_\eta z)|
 \leq \left\lfloor C^2/\eta^2\right\rfloor.
\end{equation}
Moreover,
\[
|q(u,v)-q(u',v')|
\leq\lVert u-u'\rVert_\infty+\lVert v-v'\rVert_\infty,
\]
so for all \(x,y\in H\),
\begin{equation}\label{eq:q-truncation}
|q(x+y)-q(T_\eta x+T_\eta y)|\leq4\eta.
\end{equation}
Indeed, every discarded coordinate has modulus less than \(\eta\). In each
copy of \(\ell_2\), the difference between \(x+y\) and
\(T_\eta x+T_\eta y\) therefore has sup norm at most \(2\eta\), and the
preceding Lipschitz estimate gives \eqref{eq:q-truncation}.

\begin{lemma}\label{lem:sparse-trace}
Fix nonnegative integers \(m,k\). Let \(x_1,\ldots,x_n\in H\) have support
size at most \(m\), and fix arbitrary thresholds \(c_1,\ldots,c_n\in\R\).
As \(y\) ranges over vectors of support size at most \(k\), the number of
traces
\[
\left(\mathbf1_{\{q(x_r+y)\geq c_r\}}\right)_{r=1}^n
\]
is at most
\begin{equation}\label{eq:polynomial-traces}
                         C_{m,k}(n+1)^{d_k},
\end{equation}
where \(C_{m,k}\geq1\) depends only on \(m,k\), and the finite exponent
\(d_k\geq0\) depends only on \(k\).
\end{lemma}

\begin{proof}
If \(n=0\), there is one empty trace. If \(k=0\), the only parameter is
\(y=0\). Assume \(n,k\geq1\).

To evaluate \(q(x_r+y)\) by a finite maximum, we need admissible pairs
containing zero coordinates as well as pairs of nonzero coordinates.
Put \(L=m+k\), write \(x_r=(u_r,v_r)\), and let
\[
D=\{1,\ldots,L+1\}\cup
  \bigcup_{r=1}^n
  \bigl(\operatorname{supp}u_r\cup\operatorname{supp}v_r\bigr).
\]
Thus \(|D|\leq nm+L+1\). For a parameter \(y=(u,v)\) of support size at
most \(k\), set
\[
E_y=D\cup\operatorname{supp}u\cup\operatorname{supp}v\cup\{\infty\}.
\]
Place \(\infty\) after every integer and set all coordinate values at
\(\infty\) equal to zero. It represents any integer larger than the finite
indices in \(E_y\). For every \(r\), we have
\begin{equation}\label{eq:q-finite-support}
q(x_r+y)=\max_{\substack{a,b\in E_y\\a<b}}
             |u_{r,a}+u_a+v_{r,b}+v_b|.
\end{equation}
To verify this identity, consider an admissible pair for \(z=x_r+y\).
A pair of two nonzero coordinates already has both indices in \(E_y\).
If only the first-copy coordinate is nonzero, use \(\infty\) as the second
index. If only the second-copy coordinate at index \(b\) is nonzero, the
first index of the original pair belongs to \(D\) when \(b\leq L+1\).
When \(b>L+1\), among the first \(L+1\) indices in \(D\) there is an
index where the first copy of \(z\) is zero, since \(z\) has support
size at most \(L\). In both cases \(E_y\) contains an admissible
pair giving the same value. A pair of two zero coordinates contributes zero.

Conversely, every pair on the
right of \eqref{eq:q-finite-support} is admissible, using an integer
beyond all finite indices for \(\infty\). The identity also applies when
coordinates cancel in \(x_r+y\).

Suppose \(y\) has \(\ell\leq k\) nonzero coordinates. Label them by
\(1,\ldots,\ell\) in increasing order of numerical index, with the first
copy preceding the second at a tie. By the \emph{order type} we mean the
copy of each labelled coordinate, the position of its numerical index
relative to \(D\), and the order and equality relations among these indices.
An index either equals an element of \(D\) or lies in one of its
complementary intervals in \(\N\), giving at most \(2|D|+1\) choices for
each label. The copy choices and the
order and equality relations among at most \(k\) labels contribute a
factor depending only on \(k\). Since \(|D|\leq nm+L+1\), summing over
\(\ell\leq k\) gives at most \(B_{m,k}(n+1)^k\) order types, for a
constant \(B_{m,k}\).

Fix one order type and write \(t_1,\ldots,t_\ell\) for the coordinate
values of \(y\). The order type determines which parameter coordinates
coincide with fixed support indices and which pairs in \(E_y\) satisfy
\(a<b\). All fixed nonzero coordinates have indices in \(D\). Thus, after
replacing the parameter indices by their labels,
\eqref{eq:q-finite-support} is the maximum of the absolute values of a
fixed finite list of affine functions of \(t\in\R^\ell\) for each \(r\).

For each \(r\), let \(\mathcal A_r\) be the set of sums of a coordinate
value of \(u_r\) and a coordinate value of \(v_r\), allowing zero in each
copy. Then \(|\mathcal A_r|\leq(m+1)^2\). Every affine function in the
list for \(r\) has the form \(\alpha\) or \(\alpha+w(t)\), where
\(\alpha\in\mathcal A_r\) and \(w(t)\) is one of
\[
t_i\quad(1\leq i\leq\ell),\qquad
 t_i+t_j\quad(1\leq i<j\leq\ell).
\]
Indeed, an admissible pair uses one coordinate from each copy, so it
contains at most two distinct parameter variables. There are at most
\[
                         \nu_k=k+\binom{k}{2}
\]
such linear expressions \(w\).

For each \(r\), an inequality \(|\alpha|\geq c_r\) involving a constant
affine function has a fixed truth value. For the nonconstant affine functions,
we have
\[
|\alpha+w(t)|\geq c_r
\quad\Longleftrightarrow\quad
 w(t)\geq c_r-\alpha\ \text{ or }\ w(t)\leq-c_r-\alpha.
\]
Collect all cut values in the finite set
\[
\mathcal C=\bigcup_{r=1}^n
       \{c_r-\alpha,-c_r-\alpha:\alpha\in\mathcal A_r\}.
\]
Since \(|\mathcal C|\leq2n(m+1)^2\), its points and complementary open intervals
partition \(\R\) into at most \(4n(m+1)^2+1\) pieces. On each piece,
every comparison with a cut value has a fixed truth value; singleton
pieces account for equality at a threshold.

For each \(r\), \(q(x_r+y)\geq c_r\) if and only if at least one term in
the maximum \eqref{eq:q-finite-support} is at least \(c_r\). For the fixed
order type, the piece containing \(w(t)\) for each of the at most
\(\nu_k\) expressions \(w\) determines the truth values of all the
nonconstant comparisons \(|\alpha+w(t)|\geq c_r\) and hence the threshold
trace. Counting all choices of pieces gives at most
\((4n(m+1)^2+1)^{\nu_k}\) traces for this order type. Multiplying by the
bound on the number of order types gives
\begin{align*}
\#\{\text{threshold traces}\}
 &\leq B_{m,k}(n+1)^k\bigl(4n(m+1)^2+1\bigr)^{\nu_k}\\
 &\leq C_{m,k}(n+1)^{k+\nu_k}.
\end{align*}
Thus we may take
\[
                         d_k=k+\nu_k=2k+\binom{k}{2}. \qedhere
\]
\end{proof}

\begin{theorem}
\label{thm:q-fat}
For every \(0\leq R_x,R_y<\infty\) and \(\eps>0\), the class
\(\mathcal Q_{R_x,R_y}\) has finite \(\eps\)-fat-shattering dimension.
\end{theorem}

\begin{proof}
Use the orthonormal basis \((e_j,0),(0,e_j)\) of \(H\), indexed by
\(\{1,2\}\times\N\) as above.
Apply Proposition~\ref{prop:sparse-norming-criterion} with
\(P_\eta=T_\eta\), \(a(C,\eta)=\lfloor C^2/\eta^2\rfloor\), and
\(\omega(\eta)=4\eta\). The estimates \eqref{eq:support-bound} and
\eqref{eq:q-truncation} give the approximation hypotheses, and
Lemma~\ref{lem:sparse-trace} gives the polynomial trace bound.
\end{proof}

\subsection{NIP of the norm formulas in the complete theories}

\begin{lemma}\label{lem:expansion-transfer}
Let \(M'\) be an expansion of a metric structure \(M\), and let
\(\varphi(x,y)\) be a formula in the language of \(M\). If \(\varphi\) is NIP in
\(\Th(M')\), then it is NIP in \(\Th(M)\).
\end{lemma}

\begin{proof}
Suppose that \(\varphi\) has IP in \(\Th(M)\), with thresholds \(r<s\) on
fixed sorts. Choose \(0<\delta<(s-r)/2\). By
Lemma~\ref{lem:finite-pattern-transfer}(ii), evaluation sets in \(M\)
determine boxes in every dimension between \(r+\delta\) and \(s-\delta\).
The same boxes occur in \(M'\), and Lemma~\ref{lem:finite-pattern-transfer}(i)
gives IP for \(\varphi\) in
\(\Th(M')\), a contradiction.
\end{proof}

\begin{proposition}
\label{prop:hilbert-sparse-nip}
Let \(K\) be a Hilbert space and let \(p\) be a seminorm satisfying
\(p(z)\leq C\lVert z\rVert_2\) for some \(C\geq0\). Suppose that, for every
\(0\leq R_x,R_y<\infty\) and
\(\eps>0\), the class
\[
 \{f_y:R_x\Ball_K\to\R:f_y(x)=p(x+y),\ \lVert y\rVert_2\leq R_y\}
\]
has finite \(\eps\)-fat-shattering dimension.  Then the norm formula is NIP
in the complete theory of
\[
 X_p=\left(K,N_p\right),
 \qquad
 N_p(z)=\sqrt{\lVert z\rVert_2^2+p(z)^2}.
\]
In particular, the conclusion holds if \(p\) satisfies
\(p(z)\leq C\lVert z\rVert_2\) and the hypotheses of
Proposition~\ref{prop:sparse-norming-criterion}.
\end{proposition}

\begin{proof}
Let \(X'_p\) be the expansion of \(X_p\) by
\(h(z)=\lVert z\rVert_2\) and \(p(z)\), restricted to each ball.
These predicates are bounded and uniformly continuous on each ball for
\(N_p\).

Fix \(N_p\)-ball sorts of radii \(R_x\) and \(R_y\) for \(x\) and \(y\).
Since \(h\leq N_p\), these balls are contained in the Hilbert balls
\(R_x\Ball_K\) and \(R_y\Ball_K\), respectively. The class of functions
\(x\mapsto p(x+y)\) on \(R_x\Ball_{X_p}\), with \(y\in R_y\Ball_{X_p}\),
therefore has finite fat-shattering dimension at every positive scale
by the hypothesis on the Hilbert balls. By
Lemma~\ref{lem:finite-pattern-transfer}(ii), IP for \(p(x+y)\) in
\(\Th(X'_p)\) would give finite boxes of a fixed positive width in
\(X'_p\) in every dimension, contrary to these fat-shattering bounds.
Thus \(p(x+y)\) is NIP in \(\Th(X'_p)\) on these sorts.

Let \(Z'\models\Th(X'_p)\), and let \(Z\) be its Banach reduct.
We also write \(h\) and \(p\) for their interpretations on \(Z\).
The triangle inequality for \(h\), its absolute homogeneity for rational
scalars, the parallelogram identity, the identity
\[
                  \lVert z\rVert_Z^2=h(z)^2+p(z)^2,
\]
and the inequalities
\[
        h(z)\leq\lVert z\rVert_Z\leq\sqrt{1+C^2}\,h(z)
\]
hold in \(Z'\), since they are expressed by sentences true in \(X'_p\).
In particular, \(h(z)=0\) implies \(\lVert z\rVert_Z=0\), so \(z=0\).
Continuity gives absolute homogeneity of \(h\) for real scalars. Thus \(h\)
is a norm on \(Z\).

The norm comparisons imply that \(Z\) is complete for \(h\). Indeed, an
\(h\)-Cauchy sequence is Cauchy for \(\lVert\cdot\rVert_Z\), and its limit
in \(Z\) is also its limit for \(h\). By the Jordan--von Neumann theorem,
the parallelogram identity makes \((Z,h)\) a Hilbert space
\cite[Theorem~I]{JordanVonNeumann1935}.

To test stability of \(h(x+y)\) on \(Z'\), take two sequences in fixed
\(Z\)-ball sorts. They are bounded for \(h\) and their closed linear span
in \((Z,h)\) is a separable Hilbert space, which embeds isometrically into
\(\ell_2\). The double-limit condition for \(\ell_2\) from
Section~\ref{subsec:hilbert-stability} gives stability on \(Z'\). Since
\(Z'\) was arbitrary, \(h(x+y)\) is stable in \(\Th(X'_p)\).

The two formulas \(h(x+y)\) and \(p(x+y)\) are therefore NIP in
\(\Th(X'_p)\). By Lemma~\ref{lem:nip-connectives}, so is
\[
 \sqrt{h(x+y)^2+p(x+y)^2}=N_p(x+y)
\]
in \(\Th(X'_p)\). Lemma~\ref{lem:expansion-transfer} then gives NIP
for the norm formula in \(\Th(X_p)\).
\end{proof}

\begin{theorem}\label{thm:y}
The norm formula is unstable on \(Y\) and NIP in
\(\Th(Y)\). In particular, it is NIP on \(Y\).
\end{theorem}

\begin{proof}
Instability is witnessed by \eqref{eq:half-graph}.  Since
\(q(z)\leq\sqrt2\lVert z\rVert_2\), Theorem~\ref{thm:q-fat} and
Proposition~\ref{prop:hilbert-sparse-nip}, applied with \(p=q\), give NIP
in \(\Th(Y)\), and hence on \(Y\).
\end{proof}

The fat-shattering bounds for \(q(x+y)\) also yield NIP for the norm
formula of Bollob\'as's renorming from Section~\ref{subsec:y-definition}.

\begin{corollary}[Bollob\'as's renorming]\label{cor:additive-norm}
For each \(\eps>0\), let \(Y_\eps=(H,N_\eps)\), where
\[
 N_\eps(z)=h(z)+\eps q(z).
\]
The norm formula is unstable on \(Y_\eps\) and NIP in
\(\Th(Y_\eps)\).
\end{corollary}

\begin{proof}
For \(x_i=(e_i,0)\) and \(y_j=(0,e_j)\), we have
\[
N_\eps(x_i+y_j)=
\begin{cases}
\sqrt2+2\eps,&i<j,\\
\sqrt2+\eps,&i\geq j.
\end{cases}
\]
The two iterated limits differ, so the norm formula is unstable.

Let \(Y'_\eps\) be the expansion of \(Y_\eps\) by the restrictions of
\(h\) and \(q\) to each ball sort.
The inequalities
\[
 h(z)\leq N_\eps(z)\leq(1+\sqrt2\eps)h(z)
\]
show that these predicates are bounded and uniformly continuous on each
ball and that every \(N_\eps\)-ball is contained in the corresponding
Hilbert ball. Theorem~\ref{thm:q-fat} and
Lemma~\ref{lem:finite-pattern-transfer}(ii) give NIP for \(q(x+y)\) in
\(\Th(Y'_\eps)\). The Hilbert norm identities, the equality
\(N_\eps=h+\eps q\), and the displayed norm comparisons hold in every
model of \(\Th(Y'_\eps)\). As in the proof of
Proposition~\ref{prop:hilbert-sparse-nip}, \(h\) is a complete Hilbert
norm in each such model, and \(h(x+y)\) is stable there.
Lemma~\ref{lem:nip-connectives} gives NIP for
\(N_\eps(x+y)=h(x+y)+\eps q(x+y)\) in \(\Th(Y'_\eps)\), and
Lemma~\ref{lem:expansion-transfer} gives NIP in \(\Th(Y_\eps)\).
\end{proof}

\section{NIP of the norm formula on finite \texorpdfstring{\(\ell_2\)}{l2}-sums}

\begin{lemma}\label{lem:local-nip-sum}
If the norm formula is NIP on separable Banach spaces \(E\) and \(F\), then
it is NIP on \(E\oplus_2F\).
\end{lemma}

\begin{proof}
Put \(G=E\oplus_2F\). Consider the multisorted structure consisting of the
Banach ball sorts of \(E\), \(F\), and \(G\), together with the coordinate
projections. Take a \(\kappa\)-saturated elementary extension with
\(\kappa>\max\{|E|,|F|,|G|\}\), and denote the corresponding Banach
spaces by \(E^*\), \(F^*\), and \(G^*\). Saturation passes to these
reducts, so they realize all local parameter types over \(E\), \(F\),
and \(G\), respectively, on fixed balls.
The coordinate projections remain linear, and the identity
\[
 \lVert z\rVert^2
 =\lVert\pi_Ez\rVert^2+\lVert\pi_Fz\rVert^2
\]
holds in this extension.

Fix ball radii \(R_x,R_y\) and a sequence
\(a_n=(u_n,v_n)\) in \(R_x\Ball_G\).
The coordinate projections map \(R_x\Ball_G\) into \(R_x\Ball_E\) and
\(R_x\Ball_F\), respectively.
By Lemma~\ref{lem:sequential-nip} applied to \(E\), pass to a subsequence
along which \(\lVert u_n+x\rVert\) converges for every
\(x\in R_y\Ball_{E^*}\). Applying the same lemma to \(F\) gives a further
subsequence along which \(\lVert v_n+y\rVert\) converges for every
\(y\in R_y\Ball_{F^*}\).

For every \(z\in R_y\Ball_{G^*}\), both projections have norm at most
\(R_y\), and hence
\[
 \lVert a_n+z\rVert
 =\left(
 \lVert u_n+\pi_Ez\rVert^2+
 \lVert v_n+\pi_Fz\rVert^2
 \right)^{1/2}
\]
converges along this subsequence. Since \(G^*\) realizes every local
parameter type over \(G\) on these balls,
Lemma~\ref{lem:sequential-nip} proves NIP on \(G\).
\end{proof}

\begin{corollary}\label{cor:xy}
The norm formula is unstable and NIP on \(X\oplus_2Y\). It has IP on
every ultrapower of \(X\oplus_2Y\) by a nonprincipal ultrafilter on \(\N\),
and hence in \(\Th(X\oplus_2Y)\). The space \(X\oplus_2Y\) is uniformly
convex and \(\sqrt3\)-isomorphic to Hilbert space.
\end{corollary}

\begin{proof}
By Theorem~\ref{thm:y}, the norm formula is NIP on
\(Y\). By Theorem~\ref{thm:x}, it is stable, and
hence NIP, on \(X\). Lemma~\ref{lem:local-nip-sum} therefore
implies that the norm formula is NIP on \(X\oplus_2Y\).

Embed the half-graph vectors and parameters from \eqref{eq:half-graph} in the
second summand, with zero first component.  Their norm values are unchanged,
so the formula is unstable on \(X\oplus_2Y\).

Let \(\mathcal U\) be a nonprincipal ultrafilter on \(\N\).
Lemma~\ref{lem:finite-pattern-transfer}(i) applied to the boxes in
\eqref{eq:x-cubes} gives IP on \(X^{\mathcal U}\). The isometric
inclusion \(X\hookrightarrow X\oplus_2Y\) induces an isometric inclusion
of the corresponding ultrapowers, so the same vectors and parameters
witness IP on \((X\oplus_2Y)^{\mathcal U}\).

On the underlying Hilbert sum, define the seminorm
\[
                 p(u,v)=\sqrt{p_{\oplus}(u)^2+q(v)^2}.
\]
Then \(p(u,v)\leq\sqrt2\lVert(u,v)\rVert_2\), and
\[
 \lVert(u,v)\rVert_{X\oplus_2Y}
 =\sqrt{\lVert(u,v)\rVert_2^2+p(u,v)^2}.
\]
Uniform convexity and the Hilbert-space isomorphism bound follow from
Lemma~\ref{lem:uniform-convexity}.
\end{proof}

\section{An elementary extension of \texorpdfstring{\(X\)}{X} with IP}

Khanaki's argument that copies of \(c_0\) give IP starts with an example
communicated to him by Michael Megrelishvili
\cite[Example~4.5]{KhanakiBanach}. Let \((e_i)\) be the standard unit vector
basis of \(c_0\), and put
\[
                 g_i(x)=\lVert x+e_i\rVert_\infty
                 \qquad(x\in\Ball_{c_0}).
\]
If \(E,F\subseteq\N\) are finite and disjoint, the characteristic vector
\(x=\boldsymbol{1}_F\) belongs to \(\Ball_{c_0}\) and satisfies
\[
                 g_i(x)=1\quad(i\in E),
                 \qquad
                 g_i(x)=2\quad(i\in F).
\]
Thus the sequence of vectors \((e_i)\), together with the parameters
\(\boldsymbol{1}_F\), witnesses IP for the norm formula. By James's
nondistortion theorem \cite{James1964}, every space isomorphic to \(c_0\)
contains almost isometric copies of \(c_0\)
\cite[Proposition~2.e.3]{LindenstraussTzafriri1977}. Khanaki used this fact to
deduce that the norm formula has IP on every Banach space containing an
isomorphic copy of \(c_0\) \cite[Observation~4.6]{KhanakiBanach}. He stated
that the converse fails and announced a counterexample for future work.
In subsequent correspondence, Khanaki recalled an unpublished modification
of Bollob\'as's example that replaces the order relation by a membership
relation and gives a uniformly convex norm equivalent to a Hilbert norm
\cite{KhanakiCommunication2026}. We give a construction using the
membership relation in Remark~\ref{rem:direct-ip}.

\begin{theorem}\label{thm:hilbertian-ip-space}
There is a separable elementary extension \(W\) of \(X\) such that the norm
formula has IP on \(W\), the space \(W\) is uniformly convex and
\(\sqrt3\)-isomorphic to Hilbert space, and \(W\) contains no isomorphic copy
of \(c_0\).
\end{theorem}

\begin{proof}
Let \(X'\) be the expansion of \(X\) by the restrictions of its equivalent
Hilbert norm \(h(z)=\lVert z\rVert_2\) to each ball sort. The predicate
\(h\) is bounded and \(1\)-Lipschitz on each ball sort for the norm of \(X\).
The norm identities for \(h\) and its comparison with
\(\lVert\cdot\rVert_X\) are expressible by sentences in the expanded
language, so they hold in every elementary extension of \(X'\).

Let \(\mathcal U\) be a nonprincipal ultrafilter on \(\N\).
By \eqref{eq:x-cubes}, evaluation sets for the norm formula on the
\(2\)-ball sorts of \(X'\) determine boxes in every dimension.
Lemma~\ref{lem:finite-pattern-transfer}(i)
therefore gives an ultrapower \(U=(X')^{\mathcal U}\) with vectors
\(a_i\), \(i\in\N\), and \(b_{E,F}\), for finite disjoint
\(E,F\subseteq\N\), satisfying
\[
 \lVert a_i+b_{E,F}\rVert\leq\sqrt3\quad(i\in E),
 \qquad
 \lVert a_i+b_{E,F}\rVert\geq\sqrt6\quad(i\in F).
\]
All these vectors belong to the \(2\)-ball sort of \(U\).

Identify \(X'\) with its diagonal elementary image in \(U\), and choose
a countable dense subset of each ball sort of \(X'\). The language is
countable, so downward L\"owenheim--Skolem \cite[Proposition~7.3]{BBHU08}
gives a separable elementary submodel \(W'\preceq U\) containing these
subsets and all the vectors \(a_i,b_{E,F}\). Each sort of \(W'\) is closed
in the corresponding sort of \(U\), so \(X'\subseteq W'\). Since both
are elementary submodels of \(U\), we have \(X'\preceq W'\). The Banach
reduct \(W\) is therefore a separable elementary extension of \(X\), and
the sequence \((a_i)\), with the parameters \(b_{E,F}\), witnesses IP on \(W\).

The norm axioms for \(h\), the parallelogram identity, and the inequalities
\[
                  h(z)\leq\lVert z\rVert_W\leq\sqrt3\,h(z)
\]
hold in \(W'\) by elementarity. As in the proof of
Proposition~\ref{prop:hilbert-sparse-nip}, \(h\) is a complete Hilbert norm.
The identity map from \((W,\lVert\cdot\rVert_W)\) to \((W,h)\) has
distortion at most \(\sqrt3\).

Finally, Lemma~\ref{lem:uniform-convexity} gives the midpoint estimate
\[
 \lVert x\rVert,\lVert y\rVert\leq1
 \quad\Longrightarrow\quad
 \left\lVert\frac{x+y}{2}\right\rVert^2
 \leq1-\frac1{12}\lVert x-y\rVert^2
\]
in \(X\). It is expressed by a sentence in the Banach language, so
\(X\preceq W\) implies that it holds in \(W\). Thus \(W\) is uniformly
convex. Its equivalent Hilbert norm makes \(W\) reflexive, so \(W\) contains
no isomorphic copy of \(c_0\).
\end{proof}

Theorems~\ref{thm:x}, \ref{thm:y}, and
\ref{thm:hilbertian-ip-space}, together with
Lemma~\ref{lem:uniform-convexity}, prove Theorem~\ref{thm:main}.

\begin{remark}\label{rem:direct-ip}
The membership relation also gives a uniformly convex Banach space whose
norm formula has IP. Lupini uses pairs of coordinate functionals of a
basic sequence indexed by the same membership relation to construct an
additive perturbation on an arbitrary infinite-dimensional Banach space
\cite[Theorem~4.1]{Lupini2026Norm}. Let \(K\) be a Hilbert
space with orthonormal basis consisting of the vectors \(f_i\),
\(i\in\N\), and \(e_F\), for finite \(F\subseteq\N\). Define
\[
 p(z)=\sup_{\substack{F\subseteq\N\text{ finite}\\ i\in F}}
 \bigl(|\langle z,f_i\rangle|+|\langle z,e_F\rangle|\bigr),
 \qquad N(z)=\sqrt{\lVert z\rVert_2^2+p(z)^2}.
\]
Since \(p(z)\leq\sqrt2\lVert z\rVert_2\),
Lemma~\ref{lem:uniform-convexity} shows that \(V=(K,N)\) is uniformly
convex and \(\sqrt3\)-isomorphic to Hilbert space. It is separable and
contains no isomorphic copy of \(c_0\).
As in \eqref{eq:x-cubes}, \(N(f_i+e_F)=\sqrt6\) when \(i\in F\) and
\(N(f_i+e_F)=\sqrt3\) otherwise, including when \(F=\varnothing\).
All the vectors \(f_i,e_F\) have \(N\)-norm at most \(\sqrt2\), so these
equalities witness IP for the norm formula in \(V\).
Theorem~\ref{thm:hilbertian-ip-space} gives such a space \(W\) with the
additional requirement that \(X\preceq W\).
\end{remark}

\section*{Acknowledgments}
The author thanks Karim Khanaki for helpful correspondence, for the
references to Garling and Iovino, and for bringing Lupini's preprint to
the author's attention.

\section*{AI declaration}
The author used ChatGPT Sol 6.1 to assist with the preparation and revision of this manuscript, including refinement of the counting argument underlying the fat-shattering bound in Section~\ref{subsec:fat-shattering}. The main ideas originated with the author, who has checked all mathematical statements and proofs and takes full responsibility for the mathematical content of the final manuscript.

\bibliographystyle{plainurl}
\bibliography{bibliography}
\end{document}